\magnification=\magstep1

\newcount\sec \sec=0
\input Ref.macros
\input math.macros
\forwardreferencetrue
\citationgenerationtrue
\initialeqmacro
\sectionnumberstrue

\def\mod{{\rm mod \,}}
\def\dist{{\rm dist}}

\title{Cutsets in infinite graphs}

\author{\'Ad\'am Tim\'ar}

\abstract{ We answer three
questions posed in a paper by Babson and Benjamini. They introduced a 
parameter $C_G$ for Cayley 
graphs $G$ that has significant application to percolation. For a minimal 
cutset of 
$G$ and a partition of this cutset 
into two classes, take the minimal distance between the two classes. The 
supremum of this number over all minimal cutsets and all partitions is 
$C_G$. 
We show that if it is finite for some Cayley graph of the group then it 
is finite for any (finitely generated) Cayley graph. 
Having an exponential bound for the number of minimal cutsets of size $n$ 
separating $o$ from 
infinity also turns out to be independent of the Cayley graph chosen.
We show a 1-ended example (the 
lamplighter group), where $C_G$ is infinite.
Finally, we give a new proof for a 
question of de la Harpe, 
proving that the number of $n$-element cutsets separating $o$ from 
infinity is finite unless $G$ is a finite extension of $\Z$.} 

\bottomII{Primary 05C25.  Secondary 20F65. }
{Cutsets, quasi-isometry.}
{Research partially supported by NSF Grant
DMS-0231224 and Hungarian National Foundation for Scientific Research
Grants TO34475.}

\bsection{Introduction}

In their paper \ref b.BB/, Babson and Benjamini 
introduce a parameter 
$C_G$ for an infinite Cayley graph $G$ in the following way. Let $G_+$ be 
the set of 
vertices in $G$ together with the set of ends of $G$. Given a subset $Y$ 
in $V(G)\cup E(G)$, let $$C(Y)= \sup_{Y_1\cup Y_2=Y} \{\dist(Y_1,Y_2)\}.$$
If $C(Y)\leq t$, then $Y$ is said to be $t$-{\it close}. Let $$C_G= \sup\, 
 C(\Pi),$$
where $\Pi$ ranges over all minimal cutsets between points $x,y\in G_+$.

For example, it is easy to check that for the square-grid we have 
$C_{\Z^2}=2$, while $C_G=3$ for the hexagonal grid $G$. The so called 
lamplighter group will turn out to be such that its $C_G$ is infinity 
(see Section 4).

An exponential bound on the number of minimal cutsets of size $n$ provides 
one with a nontrivial upper bound on the percolation critical
probability $p_c$ of $G$. Moreover, if $C_G$ is bounded for some one-ended
graph $G$ then the critical probability $p_u$ for having a unique infinite
cluster is smaller than 1, as shown in \ref b.BB/. There it is also shown
that for finitely  
presented groups $G$ (i.e., groups with finitely many defining relators)
$C_G$ is finite. We shall give a shorter and elementary proof
for
that in \ref s.appendix/. This established positively
Question 3 in \ref b.BS/ for
finitely presented groups, showing that a finitely group with one end has  
$p_u<1$. The original question, whether $p_u$ is smaller than 1 in 
finitely
generated groups with one end is still open, so it is natural to ask
whether
$C_G$ is finite for any Cayley graph $G$.

At the end of \ref b.BB/ there are three 
questions:

  Question 1: Does having some exponential bound on the number of minimal 
cutsets of size $n$ in a Cayley graph depend only on the group and not on 
the choice of generators?

  Question 2: Does the property ``$C_G$ is finite" for a Cayley
graph depend only on the group and not on the choice of generators?

  Question 3: Are there finitely generated groups with one end so that  
$C_G$ is not finite? 

  In  \ref s.1/ and \ref s.2/ we give positive answers to the first and 
the second 
question of \ref b.BB/ respectively. In both cases we prove the invariance 
not 
only for different Cayley graphs of the same group, but for any two 
graphs that have bounded degrees and that are quasi-isometric under a 
bijection. Similar, 
but lengthier, arguments would show that the same conclusions hold for 
any two quasi-isometric bounded degree graphs.  
 
In \ref s.3/ we show that the lamplighter group is an example of a group 
with 
one end and infinite $C_G$, hence answering Question 3 negatively.

Babson and Benjamini conjecture that for any Cayley graph the number of 
minimal cutsets separating a fixed vertex $o$ from infinity is 
exponentially bounded in the size of the cutset. They prove it for 
finitely presented groups. Problem VI.19 from \ref b.dlH/ is weaker: ``In 
a finitely-generated group which is not almost cyclic, does the size of 
spheres tend to infinity?" (A group is almost cyclic if it is a finite 
extension of $\Z$.) We give a short proof in \ref s.4/ for that the sizes 
of any 
sequence 
of  (distinct) cutsets tend to infinity. This have already been shown 
by Anna Erschler (see
the appendix of the second edition of \ref b.dlH/), using a result of
Coulhon and Saloff-Coste 
(\ref b.CS/).

For any set $X$ of vertices of a graph we denote by $N_n(X)$ the 
$n$-neighborhood of $X$, that is, the set of points at distance $\leq n$ 
from $X$ in the graph. Note that although $N_n(X)$ will be used for 
different graphs, it is the set $X$ that indicates which graph 
is understood.

In what follows, for simpler phrasing, cutset will always mean a set of 
edges whose deletion disconnects the graph. Everything remains valid if 
we use cutsets consisting of vertices. 

We do not always distinguish between vertex sets and subgraphs induced by 
them. 

A {\it bijective quasi-isometry} (or bi-Lipschitz map) {\it with 
constant} $m$ between the graphs 
$G=(V,E)$ and 
$G'=(V',E')$, 
is a bijection $f$ from $V(G)$ to $V(G')$ such that the following holds:

For any $x,y\in V(G)$, 
$$1/m\cdot \dist_G (x,y)\leq \dist_{G'} (f(x),f(y))\leq m\cdot\dist_G 
(x,y)$$.

If there is such a map between $G$ and $G'$, we shall say that they are 
{\it bijectively quasi-isometric}.
Note that different Cayley graphs of the same group are always bijectively 
quasi-isometric. 

Given a subset $X$ of $G$, $\partial X$ will denote its external   
boundary, that is, the set of vertices not in $X$ but adjacent to some
vertex of $X$. We use $\partial_i X$ for the inner boundary.

\bsection {The invariance of exponential bound}{s.1}

Given a graph $G$ and $o\in V(G)$, we say that there is an {\it
exponential
bound for the number
of minimal cutsets} separating $o$ from infinity if there is a
constant $\alpha$ such that the number of minimal cutsets of
size $n$ separating $o$
from infinity is at most $\alpha^n$. 

In this section we prove that having an exponential bound for the number 
of minimal cutsets in $o$ is 
invariant 
under bijective quasi-isometries.

We shall need the following lemma, which is basically Lemma 6 
in \ref 
b.BB/; a stronger bound is given in \ref b.Ke/. 

\procl l.benbab 
Let $G$ be a graph of degrees at most $d$ and $o$ a fixed 
vertex of it. 
The number of subsets of vertices of size $n$ that contain $o$ and induce 
a connected subgraph in 
$G$ is at most $d^{2n}$.\endprocl

\proof
Choose a spanning tree in an induced subgraph as above and define a 
``depth-first walk" in it: a walk that visits every vertex and goes 
through each edge at most twice. The walk determines the set of vertices, 
it has length $\leq 2n$, and in each step there are at most $d$ ways to 
continue 
such a walk. \Qed

\procl t.Q1 
  Let $\iota$ be a bijective quasi-isometry from a graph $G$ to $G'$ with 
constant $m$. 
Suppose that the degrees in $G$ and $G'$ are bounded by $d$. Let $o$ be 
some fixed vertex of $G$. Then there is an exponential bound for the 
number 
of minimal cutsets in $G$ separating $o$ from infinity 
if and only if there is an exponential bound in $G'$ for the 
number of cutsets separating $\iota(o)$ from 
infinity.\endprocl

\proof  Denote by $C_n$ 
and $C'_n$ the set of minimal cutsets of size $n$ separating $o$ and 
$o':=\iota(o)$ 
from infinity 
in $G$ and $G'$ respectively. Let $K_n$ ($K'_n$) be the set of connected 
subgraphs that arise as the connected component containing $o$ ($o'$) 
after 
removing an element of $C_n$ ($C'_n$) from $G$ ($G'$). We shall define a 
map $\phi$ from $K_n$ to $H_n:=\bigcup_{j=1}^{nd^{2m}} K'_j$. The map 
$\phi$ 
will 
have the property that for every $\chi \in H_n$ the set $\{ \kappa\in 
K_n: \phi(\kappa)=\chi\}$ has at most $c^n$ elements for some constant $c$ 
independent of $n$. The existence of $\phi$ shows that if $|K'_n|$ (and 
hence $|H_n|$) is at most exponential, then $|K_n|\leq |H_n| c^n$ is also 
at most 
exponential. Hence the theorem will follow. 

So, let $\kappa$ be an element in $K_n$. Define $\phi(\kappa)$ to be 
$N_m(\iota(\kappa))$. 
If
$x,y\in V(\kappa)$ are
adjacent
in $G$ then there is a path of length $\leq
m$ in $G'$ between $\iota(x)$ and $\iota(y)$. So the vertices of 
$\iota(\kappa)$ are in one component
of $\phi(\kappa)$, and since any other vertex in $\phi(\kappa)$ is in the 
same 
component of
$\phi(\kappa)$ as some vertex of $\iota(\kappa)$, we see that 
$\phi(\kappa)$ is 
connected. Now
$\phi (\kappa)$ contains $o'$. Moreover, 
for any (external) boundary vertex $y$ of $\iota(\kappa)$ in $G'$, 
$\iota^{-1}(y)$ is 
at distance at 
most $m$  from 
a boundary vertex of $\kappa$ in $G$. The set of vertices that are at 
distance $\leq m$ from
$\partial\kappa$ in $G$ have
cardinality $\leq d^m |\partial \kappa|$, so we get that  
$|\partial\iota(\kappa)|\leq 
d^m 
|\partial\kappa|$. Since $\phi(\kappa)$ is the 
$m$-neighborhood of $\iota(\kappa)$ in $G'$, the boundary of 
$\phi(\kappa)$ is in the 
$m$-neighborhood of $\partial \iota(\kappa)$ in $G'$, so 
$|\partial\phi(\kappa)|\leq d^{m} |\partial\iota(\kappa)|$. We get 
$|\partial\phi(\kappa)|\leq d^{2m} |\partial\kappa|$ 
from these two inequalities. Hence $\phi 
(\kappa)$ 
is indeed in $H_n$. 

What remains to be shown is that $|\{\kappa\in K_n : \phi (\kappa) 
=\chi\}|\leq 
c^n$ for any $\chi \in H_n$. 

Fix $\chi \in H_n$ and let $\tau:=\iota^{-1}(\chi)$. 
If $\phi (\kappa) =\chi$ for some  
$\kappa \in K_n$ then 
$\kappa\subset\tau$. Furthermore, the $m^2$-neighborhood of $\kappa$ in 
$G$ 
contains $\tau$ (by $N_m(\iota(\kappa))=\chi$ and the definition of $m$). 
Thus $\tau\setminus \kappa$ is contained in  
$N_{m^2}(\partial \kappa)$. Since $N_{m^2}(\partial \kappa)$ has $\leq 
d^{m^2} 
|\partial \kappa|=d^{m^2}n$ 
elements, $|\tau\setminus\kappa|\leq d^{m^2}n$. Now $G\setminus \kappa$ 
has only infinite
connected components by minimality of the cutsets. 
Hence we can get $\kappa$ from $\tau$ by removing a subgraph $S$ of $\tau$ 
of size $\leq d^{m^2}n$ and such that any connected component of $S$ 
contains some element of $\partial_i
\tau$. It thus 
suffices to show that there is an exponential bound on the number of such 
$S$'s.

So let $\cal{S}$ be the set of subgraphs of $\tau$ of size $\leq 
d^{m^2}n$ with the property that each component 
contains a 
vertex of $\partial_i\tau$. Any element $S$ of $\cal{S}$ can 
be described 
as follows. 

Let $S_1,\ldots,S_k$ be the components of $S$ and fix an element $r_i$ of 
$\partial_i\tau$ in $S_i$ for each $i$; let $R:=\{r_1,\ldots, 
r_k\}$. 
Now, if we first choose $R$ as a subset of $\partial_i\tau$, 
then choose the sizes of 
the 
$S_i$, and finally choose the actual 
subgraphs $S_i$ of $\tau$ of the given sizes and each incident to the 
corresponding element of $R$, then we obtain any possible 
$S\in \cal{S}$. One can choose $R$ as a subset of $\partial_i\tau$ (where 
$|\partial_i\tau|\leq d^{m^2}n$), in at most $2^{d^{m^2}n}$ ways. Once we 
have $\{r_1,\ldots, r_k\}=R$, we choose $|S_i|$ for each $i$ so that 
$\sum_{i=1}^k |S_i|\leq d^{m^2}n$ (by the definition of $\cal{S}$). This 
can be done in at most ${{d^{m^2}n}\choose {k}}\leq 2^{d^{m^2}n}$ ways. 
 Finally, 
we choose the particular $S_i$'s, knowing their sizes. By \ref l.benbab/, 
there are at most 
$d^{2(|S_1|+\ldots +|S_k|)}\leq d^{2d^{m^2}n}$ ways to do so. 

We got that there are $\leq 2^{d^{m^2}n} \cdot 2^{d^{m^2}n}\cdot 
d^{2d^{m^2}n}=(4d^2)^{d^{m^2}n}$ ways to choose $S$, thus $|\cal{S}|$ is 
exponentially bounded. This finishes the proof.\Qed

\bsection{The invariance of finiteness of $C_G$}{s.2}

In this section we show that if $C_G$ is infinite for
a graph $G$, then $C_{G'}$ is infinite 
for any 
graph $G'$ that is bijectively quasi-isometric to $G$. This means that the 
answer to Question 2 of \ref b.BB/ is 
positive.

Note that a minimal cutset $\Pi$
between two vertices $x$ and $y$ is also a minimal cutset between one of 
them and an end. Otherwise there are paths from $x$ and from $y$ to 
an 
end that do not intersect $\Pi$, and these could be used (by the 
definition 
of ends) to find a path between $x$ and $y$ that does not intersect $\Pi$, 
giving a contradiction. Hence the supremum in the definition of 
$C_G$ remains the same with
the
extra hypothesis that $y$ is in $G_+\setminus G$.

Consider a minimal cutset that separates a finite 
subgraph 
$X$ 
from 
infinity, and any boundary vertex $v$ of $X$. There is 
an 
infinite path starting from $v$ and going to infinity 
without intersecting $X$ in any point other than $v$, since 
any component of $G\setminus X$ is infinite by 
minimality.


\procl l.closure 
Suppose that a minimal cutset $S$ that separates the connected subgraph 
$X$ of $G$ from $\xi\in G_+\setminus G$ is not $l$-close. Then the set 
$S_n$ of 
edges that separate $N_n(X)$ from $\xi$ is a minimal cutset that is not 
$(l-2n)$-close.  \endprocl

\proof We may assume that $l\geq 2n$. Let $A$ and $B$ partition  
$S$ so that $dist(A,B)> l$. Then we have $N_n(A)\cap 
N_n(B)=\emptyset$. 
 The set of vertices in $N_n(X)$ that are incident to some edge in $S_n$ 
is in $N_n(A)\cup 
N_n(B)$, and it has nonempty intersection with both $N_n(A)$ and 
$N_n(B)$ because of our condition about the paths. Hence the partition 
generated on $S_n$ by 
$N_n(A)$ and $N_n(B)$ shows that it is not
$(l-2n)$-close.\Qed

\procl t.Q2 Suppose that $C_G$ is infinite for the graph $G$. Then 
$C_{G'}$ is infinite for any graph $G'$ that is bijectively 
quasi-isometric to $G$.\endprocl

\proof Let $\iota:V(G)\rightarrow V(G')$ be a bijective 
quasi-isometry 
with 
constant 
$m$. 
For each $k$, let $G_k$ be a connected subgraph in $G$ whose 
boundary is not $k$-close, and such that from any point of $\partial G_k$ 
there is a path to infinity not intersecting $G_k$.
 Such subgraphs 
exist by our assumption on $G$ and the remark about the definition of 
$C_G$. 
As in the previous section, the $m$-neighborhood $N_m(\iota(G_k))$ is 
connected. 
This 
implies that the set $S_k$ of edges that separate 
$N_m(\iota(G_k))$ from infinity is a minimal 
cutset. 
By \ref l.closure/, $S_k$ is not 
$(k/m-2m)$-close. So the $S_k$'s provide us 
with a 
sequence 
of 
minimal cutsets where the distances for certain partitions tend to 
infinity. This shows that $C_{G'}$ is infinite.\Qed

\remark The assumption that $G$ and $G'$ are bijectively quasi-isometric 
was not necessary in the last two sections. Basically similar arguments 
show that the conclusions remain true for any two quasi-isometric bounded 
degree graphs.

\bsection {The lamplighter group has infinite $C_G$}{s.3}

In this section we answer Question 3 in \ref b.BB/. 

The {\it lamplighter group} is defined as the semidirect product of $\Z$ 
with $\sum_{x\in \Z} \Z_2$. For elements $p_1,p_2\in \Z$ and 
$l_1,l_2\in \sum_{x\in \Z} \Z_2$, the product is defined as 
$$(p_1,l_1)(p_2,l_2):=(p_1+p_2,l_1\oplus S^{-m_1}l_2),$$ 
where $S$ is the left shift, $S(l)(i)=l(i+1)$ and $\oplus$ is 
componentwise addition $\mod 2$. One can think of the elements of the
lamplighter group as configurations where at each integer there is a lamp, 
either switched on or off, and there is a lamplighter standing at one of 
the integers. A possible set of generators is $\{(1,\omega), 
(0,\lambda)\}$, 
where $\omega$ stands for the sequence of all zeros, and $\lambda$ for the 
sequence of all zeros but a 1 in the 0'th position. The first generator 
corresponds to that the lamplighter 
moves one 
step to the right, and the second one to that he switches the 
lamp in his current position.  

\procl t.Q3 If $G$ is a Cayley graph of the lamplighter group then $C_G$ 
is infinite.\endprocl

The Diestel-Leader graph $DL(k,n)$ is constructed as follows. Let $T$ and 
$T'$ be a $k+1$-regular and an $n+1$-regular tree respectively, and 
suppose that they 
are rooted at infinity so that their vertices are arranged into levels 
corresponding to the integers. Do it so that a vertex of $T$ on the 
$i$'th level has $k$ 
children on the  $i+1$'th level and the parent on the $i-1$'th level; a 
vertex of $T'$ on the $i$'th level has $n$ children on the $i-1$'th level 
and the parent on the $i+1$'th level. Let the level of $v$ in $T$ 
(resp. $T'$) 
be denoted by $l_T(v)$ (resp. $l_{T'}(v)$). $DL(k,n)$ is defined on the 
vertex 
set $\{(x,x')\in V(T)\times V(T'): l_T(x)=l_{T'}(x')  \}$. 
There is an edge between $(x,x')$ 
and $(y,y')$ iff $x$ and $y$ are connected in $T$ and $x'$ and $y'$ are 
connected in $T'$.

It is well known that $DL(2,2)$ is isomorphic to a Cayley graph of the 
lamplighter group. Briefly, fix a biinfinite path $R$ in $T$ and a 
biinfinite path $R'$  in 
$T'$ 
so that these paths intersect each level in exactly one vertex. For 
each vertex in $T$ and $T'$, call the edge that goes to one of its 
children a 0-edge, 
and the other one a 1-edge. Do it so that the paths $R$ and $R'$ contain 
only 0-edges.
Now, the 
level of a 
vertex $(x,x')$ in 
$DL(2,2)$ is the position of the lamplighter, and what the lamplighter 
sees on his left (right), is just the sequence of 0's and 1's on the edges 
of the 
infinite path connecting $x$ ($x'$) to the root in infinity in $T$ 
($T'$), meaning the infinite 
path that always goes towards parents. 
 For more details, see \ref b.Wo2/.

Fix a vertex $o$ on the $0$'th level of $T$ and a vertex $o'$ on the 
$k$'th level of  $T'$. Let the subtree $F_k$ ($F'_k$) consist of the 
offspring of $o$ ($o'$) of distance at most $k$ from it. Let $H_k$ be the 
subgraph 
$\{(x,x'): x\in F_k, x'\in F'_k\}$ in $DL(2,2)$ and $C_k$ be the set of 
its 
boundary edges. Notice that $C_k$ is a minimal cutset and that it is 
the disjoint union of edges incident to
$A_k=\{(x,x'): deg_{F_k}(x)=1\}$ and
$B_k=\{(x,x'): deg_{F_k}(x')=1\}$ respectively. The distance of $A_k$ and 
$B_k$ is 
obviously $k$. 
Hence the sequence $C_k$ shows that $C_{DL(2,2)}$ is infinite.

\bsection{$C_G$ in finitely presented groups}{s.appendix}

The result in this section is the key in \ref b.BB/ to proving $p_u<1$ for 
finitely presented groups with one end. We present a shorter proof here, 
using elementary arguments and terminology, as opposed to their argument 
using cohomology groups.

We use the obvious correspondence between subsets of $E(G)$ and the
elements of $\{0,1\}^{E(G)}$ regarded as vectors, where $\mod 2$ addition
on these later corresponds to symmetric differences in the case of the
subsets.
Given a set $K$ of cycles in a graph $G$, we say that a cycle $C$ in $G$
is {\it generated by} $K$ if $C$ can be 
written as a $\mod 2$ sum of 
cycles 
from $K$. Note that 
any cycle in a Cayley graph of a finitely presented group is generated by 
the set of cycles of length at most $t$, where $t$ is the maximal length 
of relators in this presentation of the group. 

\procl t.benbab Let $G$ be a graph such that any cycle is generated by 
a set $K$ of cycles in $G$. Suppose that any cycle in $K$ has length at 
most $t$. Let $\Pi$ be a minimal cutset separating a 
vertex $x$ from $y\in V(G)\cup\{\infty\}$. Then for any nontrivial 
partition 
$\Pi_1\cup\Pi_2$ of $\Pi$ there are vertices $x_i\in \Pi_i$ ($i=1,2$) such 
that $\dist (x_1,x_2)\leq t/2$. I.e., $C_G\leq t/2$.\endprocl

\proof It is enough to show that there is a cycle in $K$ that intersects 
both $\Pi_1$ and $\Pi_2$. By minimality of the cutset $\Pi$, there are 
paths 
$P_i$ between $x$ and $y$, $i=1,2$, such that $P_i$ does not intersect 
$\Pi_{3-i}$. 
 We may write $P_1-P_2$ as a $\mod 2$ sum of cycles 
from $K$: $P_1-P_2=\sum_{c\in K'} c$, for a certain $K'\subset K$. Let 
$K'_1$ be the set of those 
cycles in $K'$ that intersect $\Pi_1$, and $K'_2:=K'\setminus K'_1$. 
Define 
$$\theta :=P_1-\sum_{c\in K'_1} c=P_2+\sum_{c\in K'_2} c.$$
The right hand side is the sum of cycles and paths that do not intersect 
$\Pi_1$, hence $\theta$ does not contain any edge of $\Pi_1$. On the 
other hand, since 
the only odd degrees that $\theta$ has are in $x$ and $y$, these two 
have to be in the same connected component of $\theta$. Thus, there is a 
path 
from $x$ to $y$ that intersects $\Pi_2$. Since $P_1$ does not, we deduce 
that some cycle in $K'_1$ does. This cycle intersects both $\Pi_1$ and 
$\Pi_2$. \Qed

\bsection{Sizes of cutsets tend to infinity}{s.4}

By the growth of a Cayley graph we mean the function that takes in $n$ the 
size 
of a ball of radius $n$.

For a rooted tree $T$ and $x\in V(T)$ let the subtree $T_x$ be defined as 
the set of all descendants of $x$ (including $x$ itself) and the edges 
induced by them.

A {\it lexicographically minimal spanning tree} $T$ of a Cayley graph $G$ 
is 
a subtree rooted at the origin and defined in 
the 
following way. Fix a linear ordering of the generators of $G$  and their 
inverses. For any vertex $v$ of the Cayley graph choose the word 
representing $v$ that is lexicographically minimal among all such words. 
There is a path in $G$ that represents this word; define $T$ to be the 
union of all these 
paths (as $v$ ranges through every vertex). The graph we get is indeed a 
spanning tree. Moreover, it is {\it subperiodic}, that is, for any $x$ in 
$V(T)$ 
there is an embedding of $T_x$ into $T$ that maps $x$ to $o$. The 
growth rates of $T$ and $G$ are the same. These are 
straightforward corollaries of the definition of $T$; for more details 
about lexicographically minimal spanning trees, see, for example, \ref 
b.LP/.

\procl l.linear If a subperiodic tree $T$ has finitely many infinite 
rays then it has linear growth.\endprocl 

\proof For $x\in V(G)$ denote by $F_x$ the union of the finite components 
of 
$T\setminus\{x\}$ not containing $o$. Define $S:=\{x\in V(T): |F_x|>|F_y|$ 
for every $y$ where $dist(o,y)<dist(o,x)\}$. By subperiodicity of $T$, 
for any $x\in V(T)$ there is an embedding $\phi$ of $T_x$ into $T$ such 
that $\phi(x)=o$. Fix $x\not = o$ and a corresponding $\phi$. If $z\in 
S\cap T_x$ then 
$\phi (F_z)$ 
can not be a 
subset of $F_{\phi(z)}$ by the definition of $S$. So $\phi(F_z)$ has the 
property that some of its vertices are mapped into vertices of an infinite 
ray starting from $\phi (z)$ and not intersecting $o$. This infinite ray 
contains no other vertices from $\phi (T_x)$ but those few from $\phi 
(F_z)$. So there are at most as many different ``$F_z$'s" as pairwise 
disjoint infinite rays, that is, $|S\cap T_x|$ is finite. This can hold 
for any $x\not = o$ iff $|S|$ is finite.
 Hence $|F_x|$ is 
bounded for every $x$, and $T$ has linear growth. \Qed

\procl t.dlH In any Cayley-graph $G$ and for $n>0$, there are only 
finitely many cutsets 
of   
size
$n$ separating a fixed vertex $o$ from infinity, unless $G$ is a
finite 
extension of $\Z$.\endprocl

\proof  Let $S_r$ stand for the sphere of
radius $r$ around $o$; $B_r$ for the ball. By a cutset 
we always 
mean a cutset separating $o$ from infinity. 

If $G$ has infinitely many ends then it is 
nonamenable and has an exponential bound for the number of minimal 
cutsets of size $n$ by \ref l.benbab/, as shown in \ref b.BB/.
So we may assume that 
$G$ has one end. Fix $n>0$.
Choose a lexicographically minimal spanning tree $T$ in $G$.
If $T$ has finitely many ends then
the group grows linearly by \ref l.linear/. Hence it is a finite 
extension of
$\Z$. (For a proof of that groups of linear growth are finite extensions 
of $\Z$, see, for example, Corollary 3.18 in \ref b.Wo1/.)
If $T$ has infinitely many ends then there is a 
ball around $o$
such that any minimal cutset of size $n$ intersects it, namely a 
ball such that there are at least $n+1$ disjoint infinite rays starting 
from its boundary. Choose $X$ to be
a set of edges that occurs in infinitely many minimal cutsets of size 
$n$ (to 
prove by contradiction), and maximal with this property. (So $0<|X|<n$.)
Since $X$ is not a cutset and $G$ has one end, there are numbers $R$ and 
$r$ 
such
that any vertex in $S_r$ is connected to $o$ by a path in $G\setminus X$ 
and not
intersecting $S_R$. But then a cutset that has no edge in $B_R$ but 
those 
of
$X$ cannot be minimal (it is necessarily a cutset without $X$ too). This
shows, by the choice of $X$, that there cannot be infinitely many minimal 
cutsets 
of
size $n$ and containing $X$. This contradiction finishes the proof. \Qed

\medbreak 
\noindent {\bf Acknowledgements.}\enspace
I thank Russell Lyons and G{\' a}bor Pete for their comments on 
the manuscript.

\startbib

\bibitem[BB]{BB} E. Babson and I. Benjamini. Cut sets and 
normed cohomology with applications to percolation. {\it Proc. Amer. 
Math. Soc.}, 127:589-597, 1999.
\bibitem[BS]{BS} I. Benjamini and O. Schramm. Percolation beyond $\Z^d$. 
{\it Electr. Commun. Prob.}, 1:71-82, 1996.

\bibitem[CS]{CS} T. Coulhon and L. Saloff-Coste. Isop\'erim\'etrie sur les 
groupes et les vari\'et\'es. {\it Rev. Mat. Iberoamericana}, 9:293-314, 
1993.

\bibitem[dlH]{dlH} P. de la Harpe. {\it Topics in Geometric Group Theory}. 
Chicago Lectures in Mathematics Series, 2000.
\bibitem[Ke]{Ke} H. Kesten. {\it Percolation Theory for Mathematicians}. 
Birkh{\"a}user, Boston.

\bibitem[LP]{LP} R. Lyons and Y. Peres, {\it Probability on Trees and 
Networks}, in preparation. 
$\;\;\;\;\;$ http://mypage.iu.edu/$\sim$rdlyons/prbtree/prbtree.html
\bibitem[Wo1]{Wo1} W. Woess. {\it Random Walks on Infinite Graphs and 
Groups}. Cambridge Tracts in Mathematics, Vol. 138, 2000.
\bibitem[Wo2]{Wo2} W. Woess. Lamplighters, Diestel-Leader graphs, random 
walks, and harmonic functions. {\it Combinatorics, Probability and
Computing}, to appear.

\endbib

\bibfile{\jobname}
\def\noop#1{\relax}
\input \jobname.bbl

\filbreak
\begingroup
\eightpoint\sc
\parindent=0pt\baselineskip=10pt

Department of Mathematics,
Indiana University,
Bloomington, IN 47405-5701
\emailwww{atimar@indiana.edu}{}
\endgroup

\bye